\documentclass[11pt]{article}
\usepackage[utf8]{inputenc}
\usepackage{amsmath,amssymb,epsfig,bbm}
\usepackage{stmaryrd,mathabx}
\usepackage{comment}
\usepackage{color}
\usepackage[T1]{fontenc}

\usepackage{interval}

\usepackage[textsize=small]{todonotes}
\usepackage{enumitem}
\usepackage{varwidth}
\setlist{nolistsep}
\usepackage{hyperref}

\newtheorem{defi}{Definition}[section]
\newtheorem{prop}[defi]{Proposition}
\newtheorem{theo}[defi]{Theorem}
\newtheorem{theofr}[defi]{Théorème}
\newtheorem{conj}[defi]{Conjecture}
\newtheorem{lemm}[defi]{Lemma}
\newtheorem{lemmfr}[defi]{Lemme}
\newtheorem{coro}[defi]{Corollary}
\newtheorem{rema}[defi]{Remark}
\newtheorem{exem}[defi]{Example}
\newtheorem{exems}[defi]{Examples}

\newcommand{\bdefi}{\begin{defi}}
\newcommand{\edefi}{\end{defi}}
\newcommand{\bprop}{\begin{prop}}
\newcommand{\eprop}{\end{prop}}
\newcommand{\btheo}{\begin{theo}}
\newcommand{\etheo}{\end{theo}}
\newcommand{\btheofr}{\begin{theofr}}
\newcommand{\etheofr}{\end{theofr}}
\newcommand{\blemm}{\begin{lemm}}
\newcommand{\elemm}{\end{lemm}}
\newcommand{\blemmfr}{\begin{lemmfr}}
\newcommand{\elemmfr}{\end{lemmfr}}
\newcommand{\brema}{\begin{rema}}
\newcommand{\erema}{\end{rema}}
\newcommand{\bexer}{\begin{exem}}
\newcommand{\eexer}{\end{exem}}
\newcommand{\bexems}{\begin{exems}}
\newcommand{\eexems}{\end{exems}}
\newcommand{\bconj}{\begin{conj}}
\newcommand{\econj}{\end{conj}}
\newcommand{\bcoro}{\begin{coro}}
\newcommand{\ecoro}{\end{coro}}
\newcommand{\dem}{\noindent{\bf Proof. }}

\usepackage{mathrsfs}
\renewcommand\mathcal{\mathscr}

\newcommand{\F}{{\cal F}}

\newcommand{\I}{{\cal I}}

\newcommand{\OOO}{{\cal O}}

\newcommand{\maths}[1]{{\mathbb #1}}  

\newcommand{\CC}{\maths{C}}

\newcommand{\HH}{\maths{H}}

\newcommand{\NN}{\maths{N}}

\newcommand{\QQ}{\maths{Q}}
\newcommand{\RR}{\maths{R}}

\newcommand{\ZZ}{\maths{Z}}

\newcommand{\aaa}{{\mathfrak a}}
\newcommand{\bbb}{{\mathfrak b}}
\newcommand{\ccc}{{\mathfrak c}}

\newcommand{\hhh}{{\mathfrak h}}

\newcommand{\mmm}{{\mathfrak m}}

\newcommand{\ppp}{{\mathfrak p}}

\newcommand{\ra}{\rightarrow}
\newcommand{\bs}{\backslash}

\newcommand{\wt}[1]{{\widetilde{#1}}}

\newcommand{\Ga}{\Gamma}

\newcommand{\cqfd}{\hfill$\Box$}

\newcommand{\bigO}{\operatorname{O}}

\newcommand{\card}{{\operatorname{Card}}}

\newcommand{\covol}{\operatorname{covol}}

\newcommand{\diam}{{\operatorname{diam}}}

\newcommand{\Nr}{\operatorname{{\tt N}}}

\renewcommand{\Re}{{\operatorname{Re}}}

\newcommand{\Vol}{\operatorname{Vol}}

\newcommand{\htr}{{\HH}^3_\RR}

\newcommand{\PSL}{\operatorname{PSL}}

\newcounter{fig}

\title{Sectorial Mertens and Mirsky formulae \\
  for imaginary quadratic number fields}
\author{Jouni Parkkonen \and Fr\'ed\'eric Paulin} 
\date{\today}

\begin{document}
\bibliographystyle{../alphanum}
\maketitle
\begin{abstract} 
We extend formulae of Mertens and Mirsky on the asymptotic behaviour
of the standard Euler function to the Euler functions of principal
rings of integers of imaginary quadratic number fields, giving versions
in angular sector and with congruences.
\end{abstract}

\section{Introduction}
\label{sec:intro}

Let $K$ be a number field of degree $n_K$, with ring of integers
$\OOO_K$, number of real places $r_1$, number of complex conjugated
places $r_2$, regulator $R_K$, class number $h_K$, number of units
$\omega_K$, discriminant $D_K$ and Dedekind zeta function $\zeta_K$
(see for instance \cite{Narkiewicz04}). Let $\I^+_K$ be the semigroup
of nonzero ideals of $\OOO_K$, let $\varphi_K :\I^+_K\ra\NN$ be the
Euler function of $K$, and let $\Nr :\I^+_K\ra\NN$ be the norm, with
$\varphi_K (a)=\varphi_K (a\OOO_K)$ and $\Nr (a)=\Nr (a\OOO_K)$ for
every $a\in \OOO_K-\{0\}$. As usual, $\ppp$ below ranges over prime
ideals in $\I^+_K$.  The functions $\bigO(\cdot)$ below depend only on
$K$.
    
Our first result (see Section \ref{sect:Mertens}) is a Mertens formula
with congruences for number fields. Though probably well-known at
least when $\mmm=\OOO_K$, we provide a proof for lack of reference
(compare with \cite[Satz 2]{Grotz79}, \cite[\S 4.3]{Cosentino99},
\cite[Theo.~3.1]{ParPau14AFST}) since arguments of its proof will
be useful for our next result. For every $\mmm\in \I_K^+$, let
$$
c_\mmm=\Nr (\mmm) \prod_{\ppp\mid\mmm} (1+\frac{1}{\Nr (\ppp)})\;.
$$
    
\btheo\label{theo:mertens} For every $\mmm\in \I_K^+$, if $n_K\geq 2$,
then as $x\ra+\infty$, we have
$$
\sum_{\aaa\in\I_K^+:\;\Nr (\aaa)\leq x,\;\mmm\mid \aaa}\varphi_K(\aaa)=
\frac{2^{r_1+r_2-1}\,\pi^{r_2}\, R_K\,h_K}
     {\omega_K\,\sqrt{|D_K|}\,\zeta_K(2)\;c_\mmm}\;x^2+
     \bigO\big(x^{2-\frac{1}{n_K}}\big)\;.
$$
\etheo

Assume in the remaining part of this introduction that $K$ is
imaginary quadratic and that $\OOO_K$ is principal.  By Dirichlet's
unit theorem, these assumptions are more or less necessary (besides
$K=\QQ$) for the following sums to be well defined and finite.

We give in Section \ref{sect:Mertenssect} a version in angular sectors
of the Mertens formula given by Theorem \ref{theo:mertens}, that will
be needed in \cite{ParPau22d}.  For all $z\in\CC^\times$,
$\theta\in\;]0,2\pi]$ and $R\geq 0$, we consider the truncated angular
    sector
\begin{equation}\label{eq:defiCzthetaR}
C(z,\theta, R)=\Big\{\rho\,e^{it}z:
t\in\;]-\frac{\theta}{2},\frac{\theta}{2}],\;
0<\rho \leq \frac{R}{|z|}\Big\}\;.
\end{equation}
It is important that the function $\bigO(\cdot)$ in the following
result is uniform in $\mmm$, $z$ and $\theta$.

\btheo \label{theo:sectcount} Assume that $K$ is imaginary quadratic
with $\OOO_K$ principal. For all $\mmm\in \I_K^+$, $z\in\CC^\times$
and $\theta\in\;]0,2\pi]$, as $x\ra+\infty$, we have
\[
\sum_{a\in\mmm\cap C(z,\theta,x)}\varphi_K(a)=
\frac{\theta}{2\,\sqrt{|D_K|}\,\zeta_K(2)\,c_\mmm}\;x^4+\bigO(x^3)\;.
\]
\etheo
        
Lastly, we give a uniform asymptotic formula for the sum in angular
sectors in $\CC$ of angle $\theta$ of the products of two shifted
Euler functions with congruences, that will be needed in
\cite{ParPau22d}. When $K=\QQ$ (the sectorial restriction is then
meaningless), this formula is due to Mirsky \cite[Thm.~9,
  Eq.~(30)]{Mirsky49} without congruences, and to Fouvry
\cite[Appendix]{ParPau22a} with congruences. For simplicity, we give a
version without congruences and without an error term in this
introduction, see Section \ref{sect:Mirskysect} Theorem
\ref{theo:mirsky} for the general statement.

\btheo \label{theo:mirskyintro} For all $z\in\CC^\times$,
$\theta\in\;]0,2\pi]$ and $k\in\OOO_K$, as $x\ra+\infty$, we have
$$
\sum_{a\in\OOO_K\cap C(z,\theta,x)}
\varphi_K(a)\,\varphi_K(a+k)\sim
\frac{\theta}{3\,\sqrt{|D_K|}}
 \prod_{\ppp}\big(1-\frac{2}{\Nr(\ppp)^2}\big)
 \prod_{\ppp\,\mid\, k\OOO_K}
 \big(1+\frac{1}{\Nr(\ppp)(\Nr(\ppp)^2-2)}\big)\;x^6\;.
$$
\etheo

Theorems \ref{theo:sectcount} and \ref{theo:mirskyintro} are used in
\cite{ParPau22d} in order to study the correlations of pairs of
complex logarithms of $\ZZ$-lattice points in the complex line at
various scalings, when the weights are defined by the Euler function,
proving the existence of pair correlation functions. We prove in
op.~cit.~that at the linear scaling, the pair correlations exhibit
level repulsion, as it sometimes occurs in statistical physics.  A
geometric application is given in op.~cit.~to the pair correlation of
the lengths of common perpendicular geodesic arcs from the maximal
Margulis cusp neighborhood to itself in the Bianchi manifolds
$\PSL_2(\OOO_K) \bs\htr$.

\medskip
\noindent{\small {\it Acknowledgements: } This research was supported
  by the French-Finnish CNRS IEA PaCap. We thank E.~Fouvry for his
  inspirational appendix in \cite{ParPau22a}.

\section{A Mertens formula with congruences for number fields}
\label{sect:Mertens}

Recall that $\I^+_K$ is the semigroup of nonzero (integral) ideals of
the Dedekind ring $\OOO_K$ (with unit $\OOO_K$). For all $I,J\in
\I^+_K$, we write $J\mid I$ if $I\subset J$, we denote by $(I,J)=I+J$
the greatest common ideal divisor of $I$ and $J$, by $[I,J]=I\cap J$
the least common ideal multiple of $I$ and $J$, and by $IJ$ the
product ideal of $I$ and $J$.

We denote by $\Nr (I)=\card (\OOO_K/I)$ the (absolute) {\it norm} of
$I\in \I^+_K$, which is completely multiplicative. The {\it norm} of $a\in
\OOO_K-\{0\}$ is
\[
\Nr (a)=\Nr (a\OOO_K)\;.
\]
It coincides with the (relative) norm $N_{K/\QQ}(a)$ of $a$ (see for
instance \cite{Narkiewicz04}), and in particular is equal to $|a|^2$
if $K$ is imaginary quadratic.

Recall that the %(complex-analytic) 
{\it Dedekind zeta function}
$\zeta_K:\{s\in\CC:\Re(s)>1\}\ra \CC$ of $K$ is defined (see for
instance \cite[\S 7.1]{Narkiewicz04}) equivalently by
$$
\zeta_K  (s)=\sum_{\aaa\in\I^+_K}\frac{1}{\Nr (\aaa)^s}=
\prod_{\ppp}\big(1-\frac{1}{\Nr (\ppp)^s}\big)^{-1}\;.
$$

We denote by $\varphi_K:\I^+_K\ra \NN$ the {\it Euler function} of $K$,
defined (see for instance \cite[page 13]{Narkiewicz04}) equivalently
by
\[
\forall\;\aaa\in\I_K^+,\;\;\;
\varphi_K  (\aaa)=\card\big((\OOO_K/\aaa)^\times\big)=
\Nr (\aaa)\prod_{\ppp\mid \aaa}\big(1-\frac{1}{\Nr (\ppp)}\big)\;.
\]
For every $a\in \OOO_K-\{0\}$, we define $\varphi_K(a)=
\varphi_K(a\OOO_K)$. Note that the Euler function $\varphi_K$ is
multiplicative\footnote{Recall that a function $f:\I^+_K\ra\CC^\times$
is {\it multiplicative} if $f(\OOO_K)=1$ and if for all coprime
integral ideals $\aaa,\bbb$ in $\I^+_K$, we have $f(\aaa\bbb)=
f(\aaa)f(\bbb)$.}  by the Chinese remainder theorem. We have
\begin{equation}\label{eq:normEUler}
  \Nr (\aaa)=\sum_{\bbb\mid \aaa}\;\varphi_K(\bbb)\;,
\end{equation}
as checked by telescopic sum when $\aaa$ is a power of a
prime ideal, and by multiplicativity.

We denote by $\mu_K:\I^+_K\ra \ZZ$ the {\it Möbius function} of $K$,
defined  by
$$
\forall\;\aaa\in\I_K^+,\;\;\;
\mu_K(\aaa)=\begin{cases} 1& {\rm if~}\aaa=\OOO_K\\
0 & {\rm if~}\ppp^2\mid \aaa
   {\rm ~for~some~prime~ideal~}\ppp\\
   (-1)^m &{\rm if~} \aaa=\ppp_1\dots\ppp_m
   {\rm ~for~pairwise~distinct~prime~ideals~}\\&\ppp_1,\dots,\ppp_m{\rm ~and~}
   m\in\NN-\{0\}\;.
\end{cases}
$$
%The {\it Dirichlet convolution} $f*g:\I^+_K\ra \CC$ of two functions
%$f,g:\I^+_K\ra \CC$ (see for instance loc.~cit.), defined by
%$$
%f*g(\aaa)=\sum_{\bbb\mid \aaa}\;f(\bbb)\;g(\aaa\bbb^{-1})\;,
%$$
%is commutative and associative, and
For every $a\in \OOO_K-\{0\}$, we define $\mu_K(a)= \mu_K(a\OOO_K)$.
We have (see for instance \cite{Shapiro59}) %(see for instance op.~cit.) 
the {\it Möbius inversion
  formula}: for all $f,g:\I^+_K\ra \CC$, 
\begin{equation}\label{eq:Mobinvform}
  f(\aaa)=\sum_{\bbb\mid \aaa}\;g(\bbb) {\rm~~~if~and~only~if~~~}
  g(\aaa)=\sum_{\bbb\mid \aaa}\;\mu_K(\bbb)f(\aaa\bbb^{-1})\;.
\end{equation}
In particular, since the norm is completely multiplicative and by
Equation \eqref{eq:normEUler}, we have
\begin{equation}\label{eq:MobiusinvEuler}
\forall\;\aaa\in\I_K^+,\;\;\; \frac{\varphi_K(\aaa)}{\Nr (\aaa)}=
\sum_{\bbb\mid \aaa}\;\frac{\mu_K(\bbb)}{\Nr (\bbb)}\;.
\end{equation}

\medskip
\noindent{\bf Proof of Theorem \ref{theo:mertens}. }
In this proof, all functions $\bigO(\cdot)$ depend only on $K$.  Let
\begin{equation}\label{eq:defirhoK}
\rho_K=
\frac{2^{r_1}\,(2\pi)^{r_2}\, R_K\,h_K}{\omega_K\,\sqrt{|D_K|}}\;.
\end{equation}
Recall (see for instance \cite[Theo.~5]{MurvOrd97}) that, as
$x\ra+ \infty$, we have
\begin{equation}\label{eq:asympcard}
\card\{\aaa\in\I_K^+:\Nr (\aaa)\leq x\} =
\rho_K\,x+\bigO(x^{1-\frac{1}{n_K}})\;.
\end{equation}
By Abel's summation formula, as $y\ra+ \infty$, we have
\begin{equation}\label{eq:sumNaaaleqy}
\sum_{\aaa\,\in\I_K^+:\;\Nr (\aaa)\leq y}\Nr (\aaa)=
\sum_{1\leq n\leq y}n\,\card\{\aaa\in\I_K^+:\Nr (\aaa)=n\} =
\frac{\rho_K}{2}\,y^2+\bigO(y^{2-\frac{1}{n_K}})\;.
\end{equation}
Furthermore, we have
\[
\card\{\aaa\in\I_K^+:\Nr (\aaa)= y\} =
\bigO(y^{1-\frac{1}{n_K}})\;.
\]
This formula implies since $\Nr ((\bbb,\mmm)) \leq \Nr (\mmm)$ that
\begin{equation}\label{eq:controlrestzeta2}
  \Big|\;\sum_{\bbb\in\I_K^+:\;\Nr (\bbb)\geq x}\mu_K(\bbb)
  \frac{\Nr ((b,\mmm))}{\Nr (\bbb)^2}\;\Big|\;=
  \bigO\Big(\Nr (\mmm)\sum_{n\geq x}
  \frac{n^{1-\frac{1}{n_K}}}{n^2}\Big)
  =\bigO( \Nr (\mmm)\;x^{-\frac{1}{n_K}})\;.
\end{equation}
Let us denote by $S_\mmm(x)$ the sum on the left hand side in the
statement of Theorem \ref{theo:mertens}. Note that by the Gauss
lemma, for all $\mmm,\bbb,\ccc\in\I_K^+$, we have $\mmm\mid\bbb\ccc$
if and only if $\mmm(\mmm,\bbb)^{-1}\mid\ccc$. Then by Equation
\eqref{eq:MobiusinvEuler}, by the change of variable $\ccc=
\mmm(\mmm,\bbb)^{-1}\aaa$, by the complete multiplicativity of the
norm, by Equation \eqref{eq:sumNaaaleqy} with $y=\frac{\Nr
  ((\bbb,\mmm))x} {\Nr (\bbb) \Nr (\mmm)}$, since $\Nr ((\bbb,\mmm))
\leq \Nr (\mmm)$, and by Equation \eqref{eq:controlrestzeta2}, we have
\begin{align}
S_\mmm(x)&=\sum_{\aaa\in\I_K^+:\;\Nr (\aaa)\leq x,\;\mmm\mid \aaa}
\;\;\;\sum_{\bbb,\ccc\,\in\I_K^+:\;\bbb\ccc=\aaa}\mu_K(\bbb)\;
\Nr (\ccc)\nonumber
\\ &=\sum_{\bbb\in\I_K^+:\;\Nr (\bbb)\leq x}\mu_K(\bbb)\;\;\;
\sum_{\ccc\,\in\I_K^+:\;\Nr (\ccc)\leq
  \frac{x}{\Nr (\bbb)},\;\mmm\mid\bbb\ccc}\;
\Nr (\ccc)\nonumber
\\ &=\sum_{\bbb\in\I_K^+:\;\Nr (\bbb)\leq x}\mu_K(\bbb)\;\;\;
\sum_{\aaa\,\in\I_K^+:\;\Nr (\aaa)\leq
  \frac{\Nr ((\bbb,\mmm))x}{\Nr (\bbb)\Nr (\mmm)}}\;
\frac{\Nr(\mmm)}{\Nr ((\bbb,\mmm))}\;\Nr(\aaa)\nonumber
\\ &=\Big(\sum_{\bbb\in\I_K^+:\;\Nr (\bbb)\leq x}\mu_K(\bbb)
\frac{\Nr ((\bbb,\mmm))}{\Nr (\bbb)^2}
\;\Big) \frac{\rho_K}{2\,\Nr (\mmm)}
\;x^2+\bigO(\,x^{2-\frac{1}{n_K}})\nonumber
\\ &=\Big(\sum_{\bbb\in\I_K^+}\mu_K(\bbb)
\frac{\Nr ((\bbb,\mmm))}{\Nr (\bbb)^2}\;
\Big) \frac{\rho_K}{2\,\Nr (\mmm)}
\;x^2+\bigO(\,x^{2-\frac{1}{n_K}})
\;.\label{eq:premtransfSmx}
\end{align}
By decomposing a nonzero integral ideal $\bbb$ into powers of prime
ideals, by the definition of the Möbius function, and by the Euler
product formula for the Dedekind zeta function, we have
\begin{align*}
\sum_{\bbb\in\I_K^+}\mu_K(\bbb)
\frac{\Nr ((\bbb,\mmm))}{\Nr (\bbb)^2}&=
\prod_{\ppp\,\nmid\, \mmm}(1-\frac{1}{\Nr (\ppp)^2})
\prod_{\ppp\,\mid\,\mmm}(1-\frac{1}{\Nr (\ppp)})=
\frac{1}{\zeta_K(2)}\prod_{\ppp\,\mid\,\mmm}
\frac{\Nr (\ppp)}{1+\Nr (\ppp)}\;.
\end{align*}
Equations \eqref{eq:premtransfSmx} and \eqref{eq:defirhoK} hence imply
Theorem \ref{theo:mertens}.
\cqfd

\section{A sectorial Mertens formula}
\label{sect:Mertenssect}

Assume in the remaining part of this paper that $K$ is imaginary quadratic
and that $\OOO_K$ is principal (or equivalently factorial (UFD)). By
Dirichlet's unit theorem, the group of units $\OOO_K^\times$, whose
order we denote by $|\OOO_K^\times|$, is finite if and only if
$(r_1,r_2)$ is equal to $(1,0)$ or $(0,1)$. This justifies our
restriction, the case $K=\QQ$ being well-known. With the notation of
the beginning of the introduction, we then have (see for instance
\cite{Narkiewicz04}) $D_K\in\{-4,-8,-3,-7,-11,-19,-43,-67,-163\}$, and
\begin{equation}\label{eq:constKimagOKprinc}
r_1=0, \quad r_2=1, \quad n_K=2, \quad R_K=1, \quad
\omega_K=|\OOO_K^\times| \quad{\rm and}\quad h_K=1\;.
\end{equation}

Given a {\it $\ZZ$-lattice} $\vec\Lambda$ in the Euclidean space $\CC$
(that is, a discrete (free abelian) subgroup of $(\CC,+)$), we denote
by $\covol_{\vec\Lambda}=\Vol(\CC/ \vec\Lambda)$ the area of a
fundamental parallelogram $\F_{\vec\Lambda}$ for $\vec\Lambda$ and by
$\diam_{\vec\Lambda}$ the diameter of $\F_{\vec\Lambda}$.  Note that
every element $\mmm\in \I^+_K$ is a $\ZZ$-lattice in $\CC$ with
\begin{equation}\label{eq:calccocoldiamideal}
\covol_\mmm=\Nr (\mmm)\covol_{\OOO_K}=\frac{\Nr (\mmm)\,\sqrt{|D_K|}}{2}
\;\;\;{\rm and}\;\;\;
\diam_\mmm=\bigO(\sqrt{|D_K|\,\Nr (\mmm)}\;)
\end{equation}
since $\diam_{\OOO_K}=|1+\frac{\sqrt{D_K}}{2}|$ if $D_K\equiv
0\!\!\mod 4$ and $\diam_{\OOO_K}=|\frac{3+\sqrt{D_K}}{2}|$ if
$D_K\equiv 1\!\!\mod 4$.

With the notation of Equation \eqref{eq:defiCzthetaR}, note that for
every $z'\in\CC^\times$, we have
\begin{equation}\label{eq:dilatCzthetaR}
  z'C(z,\theta,R)=C(zz',\theta, R\,|z'|\,)\;.
\end{equation}

\medskip
\noindent{\bf Proof of Theorem \ref{theo:sectcount}. }
Let $z\in\CC^\times$, $\theta\in\;]0,2\pi]$ and $y>0$. Since
$\operatorname{Area} (C(z,\theta,y))= \frac{\theta}{2}\,y^2$, the
standard Gauss counting argument, the finiteness of the
number of  imaginary quadratic number fields with class number $1$, and
the equality on the left of Formula \eqref{eq:calccocoldiamideal} give
\begin{align*}
\card\big(\OOO_K\cap C(z,\theta,y)\big)&=
\frac{\operatorname{Area}(C(z,\theta,y))}{\covol_{\OOO_K}}+
\bigO\big(\frac{\diam_{\OOO_K}\;y}{\covol_{\OOO_K}}\big)\\ &=
\frac{\theta}{\sqrt{|D_K|}}\,y^2+
\bigO(y)\;.
\end{align*}
Since the map $z'\mapsto |z'|^2=\Nr (z')$ takes only
integral values on $\OOO_K$, by Abel's summation formula, as
$y\ra+ \infty$, we have
\begin{align}
\sum_{d\in\OOO_K\cap C(z,\theta,y)}|d|^2&=
\sum_{1\leq n\leq y^2}n\,\card\{d\in\OOO_K\cap C(z,\theta,y):|d|^2=n\}
\nonumber\\&=
\frac{\theta}{2\,\sqrt{|D_K|}}\,y^4+
\bigO(y^3)\;.
\label{eq:sumNaaaleqysect}
\end{align}

For all $x\geq 1$ and $\bbb\in\I^+_K$, let us fix $b, m ,(b, m ) \in
\OOO_K-\{0\}$ such that $\bbb=b\OOO_K$, $\mmm= m \OOO_K$ and
$(\bbb,\mmm)=(b, m )\OOO_K$. Since for every $c\in\OOO_K-\{0\}$ we
have $ m \mid bc$ if and only if $\frac{ m }{ (b, m )}\mid c$, by the
change of variable $c=\frac{ m }{ (b, m )}\,d$, by Equation
\eqref{eq:dilatCzthetaR} and by Equation \eqref{eq:sumNaaaleqysect}
applied with $y=\frac{x|(b, m )|} {| m |\,|b|}$, if
\[
S_\bbb= \sum_{\ccc\in \I^+_K,\,a\in\mmm\cap C(z,\theta,x)\;:\;
\bbb\ccc=a\OOO_K}\Nr (\ccc)\;,
\]
we have
\begin{align*}
  S_\bbb&= \sum_{c\,\in\,\OOO_K-\{0\},\,a\in\mmm\cap
    C(z,\theta,x)\;:\; bc=a} |c|^2\\ &=
  \sum_{c\,\in\,\OOO_K-\{0\}\;:\;bc\,\in \,C(z,\theta,x),\; m \,\mid\,
    bc} |c|^2 =\sum_{d\,\in\,\OOO_K-\{0\}\;: \;d\,\in \,C(\frac{z(b, m
      )}{ m \,b},\,\theta, \frac{x|(b, m )|}{| m |\,|b|}),\;}
  \frac{\Nr (\mmm)}{\Nr ((\bbb,\mmm))} \;|d|^2\\ &=\frac{\theta\,\Nr
    ((\bbb,\mmm))} {2\,\sqrt{|D_K|}\,\Nr (\mmm)\, \Nr (\bbb)^2} \;x^4
  +\bigO\Big(\frac{x^3}{\Nr (\bbb)^{3/2}}\Big)\;.
\end{align*}
Let us denote by $S_{\mmm,z,\theta}(x)$ the sum on the left hand side
in the statement of Theorem \ref{theo:sectcount}. Then by Equation
\eqref{eq:MobiusinvEuler}, we have
\begin{align*}
S_{\mmm,z,\theta}(x)&= \sum_{a\in\mmm\cap C(z,\theta,x)}\varphi_K(a\OOO_K)
  =\sum_{a\in\mmm\cap C(z,\theta,x)}\;\;\;
  \sum_{\bbb,\ccc\,\in\I_K^+:\;\bbb\ccc=a\OOO_K}\mu_K(\bbb)
  \Nr (\ccc)\\&=
  \sum_{\bbb\in\I_K^+:\;\Nr (\bbb)\leq x^2}\mu_K(\bbb)\;\;S_\bbb
\\ &=\Big(\sum_{\bbb\in\I_K^+:\;\Nr (\bbb)\leq x^2}\mu_K(\bbb)\;
\frac{\Nr ((\bbb,\mmm))}{\Nr (\bbb)^2}\;
\Big) \frac{\theta}{2\,\sqrt{|D_K|}\,\Nr (\mmm)}
\;x^4+\bigO(x^3)\;.
\end{align*}

The proof then proceeds exactly as in the proof of Theorem
\ref{theo:mertens}.
\cqfd

\section{A sectorial Mirsky formula}
\label{sect:Mirskysect}

We now give a uniform asymptotic formula for the sum in angular
sectors of the products of shifted Euler functions with
congruences. For all $z\in\CC^\times$, $\theta\in\;]0,2\pi]$,
$k\in\OOO_K$, $\mmm\in \I_K^+$ and $x\geq 1$, let
\begin{equation}\label{eq:defiSzthetakmmmx}
S_{z,\theta,k,\mmm}(x)=
\sum_{a\in\mmm\cap C(z,\theta,x)}\varphi_K(a)\,\varphi_K(a+k)\;.
\end{equation}

\btheo \label{theo:mirsky} Assume that $K$ is imaginary quadratic with
$\OOO_K$ principal. There exists a universal constant $C>0$ such that
for all $k\in\OOO_K$ and $\mmm\in \I_K^+$, there exists
$c_{\mmm,k}\in\;]0,1]$ such that for all $z\in\CC^\times$,
    $\theta\in\;]0,2\pi]$ and $x\geq 1$, we have
$$
\Big|\;S_{z,\theta,k,\mmm}(x)-
\frac{\theta\,c_{\mmm,k}}{3\,\sqrt{|D_K|}}\;x^6\;\Big|\leq
C\big((1+\sqrt{\Nr(k)}\,)\,x^5+\Nr(k)\,x^4\big)\;.
$$
\etheo

\medskip
We will prove Theorem \ref{theo:mirsky} at the end of this Section
after giving a number of Lemmas required for the proof.  We fix
$k\in\OOO_K$ and $\mmm=m\OOO_K\in \I_K^+$, and we define
$\hhh=k\OOO_K$, which is a possibly zero integral ideal. We start by
giving the first definition and a simpler formula for the constant
$c_{\mmm,k}$ that appears in the statement of Theorem
\ref{theo:mirsky}. We define
\begin{equation}\label{eq:deficmk}
  c_{\mmm,k}=\sum_{\substack{\bbb,\ccc\in \I_K^+\\
      (\bbb,\ccc)\,\mid\, \hhh,\;
      (\ccc(\bbb,\mmm), \mmm(\bbb,\ccc))\, \mid\, \hhh\bbb}}
  \mu_K(\bbb)\,\mu_K(\ccc)\;\frac{\Nr\big((\ccc(\bbb,\mmm),\,
  \mmm(\bbb,\ccc))\big)}{\Nr(\bbb)^2\,\Nr(\ccc)^2\,\Nr(\mmm)}\;,
\end{equation}
and
\[
c'_{\mmm}=\inf_{k\in\OOO_K}c_{\mmm,k}\;.
\]

\blemm\label{lem:controlconstantcmmmk} The series in Equation
\eqref{eq:deficmk} defining $c_{\mmm,k}$ converges absolutely. We have
$c_{\mmm,k}\leq 1$ and $c'_{\mmm}>0$.  Furthermore, we have
\begin{equation}\label{eq:simplifcmk}
c_{\mmm,k}= \frac{1}{\Nr(\mmm)}\prod_{\substack{\ppp\\(\ppp,\mmm)\,\mid\, \hhh}}
\big(1-\frac{\Nr((\ppp,\mmm))}{\Nr(\ppp)^2}\big)\;
\prod_{\ppp} \big(1-\frac{\kappa_{\mmm,\hhh}(\ppp)\;
  \kappa'_\hhh(\ppp)\,\Nr((\ppp,\mmm))}{\Nr(\ppp)^2}\big)\;,
\end{equation}
where
\begin{equation}\label{eq:defikappa}
\kappa_{\mmm,\hhh}(\ppp)=\Big\{\begin{array}{l}
(1-\frac{\Nr((\ppp,\mmm))}{\Nr(\ppp)^2})^{-1}{\rm ~~~if~~~}
(\ppp,\mmm)\mid \hhh\\ 1
{\rm ~~~otherwise} \end{array} {\rm ~~~and~~~}
\kappa'_{\hhh}(\ppp)=\left\{\begin{array}{l}
1-\frac{1}{\Nr(\ppp)}{\rm ~~~if~~~} \ppp\mid \hhh\\ 1 {\rm ~~~otherwise.}
\end{array}\right.
\end{equation}
\elemm

In the special case $\mmm=\OOO_K$, Equation \eqref{eq:simplifcmk}
becomes
{\small\begin{align}
  c_{\OOO_K,k}&=\prod_{\ppp}\big(1-\frac{1}{\Nr(\ppp)^2}\big)
  \prod_{\ppp\,\mid\, \hhh}\Big(1-\frac{(1-\frac{1}{\Nr(\ppp)^2})^{-1}
    (1-\frac{1}{\Nr(\ppp)})}{\Nr(\ppp)^2}\Big) \prod_{\ppp\,\nmid\,
    \hhh}\Big(1-\frac{(1-\frac{1}{\Nr(\ppp)^2})^{-1}}
       {\Nr(\ppp)^2}\Big)\nonumber\\ & =
       \prod_{\ppp}\big(1-\frac{1}{\Nr(\ppp)^2}\big)
       \prod_{\ppp\,\mid\,
         \hhh}\Big(1-\frac{\Nr(\ppp)-1}{\Nr(\ppp)(\Nr(\ppp)^2-1)}\Big)
       \prod_{\ppp}\big(1-\frac{1}{\Nr(\ppp)^2-1}\big)
       \prod_{\ppp\,\mid\,
         \hhh}\Big(1-\frac{1}{\Nr(\ppp)^2-1}\Big)^{-1}\nonumber \\ &
       =\prod_{\ppp}\big(1-\frac{2}{\Nr(\ppp)^2}\big)
       \prod_{\ppp\,\mid\,
         \hhh}\big(1+\frac{1}{\Nr(\ppp)(\Nr(\ppp)^2-2)}\big)\;.
\label{eq:computCOKk}
\end{align}}

\noindent Theorem \ref{theo:mirskyintro} in the introduction follows from
Theorem \ref{theo:mirsky} and the above computation.

\medskip
\dem Let us prove that uniformly in $x\geq 1$, we have
\begin{equation}\label{eq:Zmmmk}
\sum_{\substack{\bbb,\ccc\in \I_K^+\;:\; \Nr(\bbb)\geq x,\\
      (\bbb,\ccc)\,\mid\, \hhh,\;
      (\ccc(\bbb,\mmm), \mmm(\bbb,\ccc))\, \mid\, \hhh\bbb}}
  \frac{\Nr\!\big((\ccc(\bbb,\mmm),\,
    \mmm(\bbb,\ccc))\big)}{\Nr(\bbb)^2\,\Nr(\ccc)^2\,\Nr(\mmm)}=
  \bigO\big(\frac{1}{\sqrt{x}}\big)\;.
\end{equation}
This implies, by taking $x=1$, that the first claim of Lemma
\ref{lem:controlconstantcmmmk} is satisfied, since the Möbius function
has values in $\{0,\pm1\}$. Let us denote by $Z_{\mmm,\hhh}(x)$ the above
sum. Since $\Nr\!\big((\ccc(\bbb,\mmm),\, \mmm(\bbb,\ccc))\big)\leq
\Nr(\mmm(\bbb,\ccc))$, we have
\begin{align*}
Z_{\mmm,\hhh}(x)&\leq \sum_{\bbb,\ccc\in \I_K^+\;:\; \Nr(\bbb)\geq x}
\frac{\Nr((\bbb,\ccc))}{\Nr(\bbb)^2\,\Nr(\ccc)^2}
\leq\sum_{\substack{\aaa,\bbb',\ccc'\in \I_K^+\\ \Nr(\bbb')\geq x/\Nr(\aaa)}}
\frac{\Nr(\aaa)}{\Nr(\aaa\bbb')^2\,\Nr(\aaa\ccc')^2}\\ & =
\sum_{\ccc'\in \I_K^+}\frac{1}{\Nr(\ccc')^2}\sum_{\aaa\in \I_K^+}
\frac{1}{\Nr(\aaa)^{5/2}}
\sum_{\substack{\bbb'\in \I_K^+\\\Nr(\aaa)\Nr(\bbb')\geq x}}
\frac{1}{\Nr(\bbb')^{3/2}\;(\Nr(\aaa)\Nr(\bbb'))^{1/2}}
\leq \zeta_K(2)\,
\zeta_K(\frac{5}{2})\,\zeta_K(\frac{3}{2})\;\frac{1}{\sqrt{x}}\;.
\end{align*}
Equation \eqref{eq:Zmmmk} follows, since there are only finitely many
fields $K$ satisfying the assumptions of Theorem \ref{theo:mirsky}.

\medskip
The proof of Equation \eqref{eq:simplifcmk} that we now give is
similar to Fouvry's proof of Equation (21) in
\cite[Appendix]{ParPau22a}.

For every $\bbb\in \I_K^+$, let $\chi_\bbb:\I_K^+\ra\{0,1\}$ be the
characteristic function of the set of elements $\ccc\in \I_K^+$ such
that $(\ccc,\bbb) \mid \hhh$. Let us define a map $\psi_\bbb:\I_K^+\ra
\I_K^+$ by
\begin{equation}\label{eq:defpsisubd}
  \psi_\bbb: \ccc\mapsto
  \big(\ccc, \frac{\mmm}{(\bbb,\mmm)}\,(\bbb,\ccc)\big)\;.
\end{equation}

Note that the assertion $(\ccc(\bbb,\mmm), \mmm(\bbb,\ccc)) \mid
\bbb\,\hhh$ is equivalent to the assertion
\[
\psi_\bbb(\ccc)\, \mid\, \frac{\bbb}{(\bbb,\mmm)}\,\hhh\;.
\]
For every $\bbb\in \I_K^+$, let $\chi^*_\bbb:\I_K^+\ra\{0,1\}$ be the
characteristic function of the set of elements $\ccc\in \I_K^+$ such
that the above divisibility assertion is satisfied.  Let us finally
define a map $C^*:\I_K^+\ra\RR$ (which depends on $\mmm$ and $\hhh$)
by
\begin{equation}\label{eq:defceuler}
C^*:\bbb\mapsto\sum_{\ccc\in \I_K^+}\frac{\mu_K(\ccc)}{\Nr(\ccc)^2}
\;\chi_\bbb(\ccc)\;\chi^*_\bbb(\ccc)\,\Nr(\psi_\bbb(\ccc))\;.
\end{equation}
By the absolute convergence property, Equation \eqref{eq:deficmk} then
becomes
\begin{equation}\label{eq:cabktech2}
c_{\mmm,k}=\;\frac{1}{\Nr(\mmm)}\sum_{\bbb\in \I_K^+}
\frac{\mu_K(\bbb)}{\Nr(\bbb)^2}\; \Nr((\bbb,\mmm))\;C^*(\bbb)\;.
\end{equation}
In order to transform the series $C^*(\bbb)$ defined by Formula
\eqref{eq:defceuler} into an Eulerian product and in order to analyse
it, we will use the following two lemmas.

\blemm\label{lem:multiplicatives}
For every $\bbb\in \I_K^+$, the maps $\chi_\bbb$, $\chi^*_\bbb$ and
$\psi_\bbb$ on $\I_K^+$ are multiplicative.
\elemm

\dem We have $\psi_\bbb(\OOO_K)=\OOO_K$ and
$\chi_\bbb(\OOO_K)=\chi^*_\bbb(\OOO_K)=1$. Let $I,J\in \I_K^+$ be
coprime.

The equality $(IJ,\bbb)=(I,\bbb)(J,\bbb)$ and the fact that $(I,\bbb)$
and $(J,\bbb)$ are coprime imply that $\chi_\bbb(IJ)= \chi_\bbb(I)
\chi_\bbb(J)$.

In order to prove the multiplicativity of the map $\psi_\bbb$, we
write
$$
\psi_\bbb(I J) =
\big(I J , \frac{\mmm}{(\bbb,\mmm)}\,(\bbb,I J)\big)=
\big(I, \frac{\mmm}{(\bbb,\mmm)}\,(I,\bbb)(J,\bbb)\big)
\big(J, \frac{\mmm}{(\bbb,\mmm)}\,(I,\bbb)(J,\bbb)\big)\;.
$$
Since $I$ is coprime to $(J,\bbb)$ and since $J$
is coprime to  $(I,\bbb)$, we obtain as wanted the equality
$\psi_\bbb(I J) = \psi_\bbb(I)\,\psi_\bbb(J)$.

Finally, the multiplicativity property
$\chi^*_\bbb(IJ)=\chi^*_\bbb(I)\chi^*_\bbb(J)$ of the function
$\chi^*_\bbb$ is a consequence of the multiplicativity of the map
$\psi_\bbb$ and of the fact that $\psi_\bbb(I)$ and $\psi_\bbb(J)$ are
coprime.  \cqfd

\blemm\label{lem:valeurschichippsi}
For every prime ideal $\ppp$ and every $\bbb\in \I_K^+$, we have
$$
\psi_\bbb(\ppp)= \left\{\begin{array}{l} \ppp {\rm ~~~ if~~~} \ppp\mid \bbb,\\
(\ppp,\mmm){\rm ~~~ otherwise},\end{array}\right.
$$
and
$$
\chi_\bbb(\ppp)\;\chi^*_\bbb(\ppp)=1\Leftrightarrow \left\{\begin{array}{l}
\ppp\mid(\bbb,\hhh) %{\rm ~~~ and~~~} \ppp\mid\frac{\bbb}{(\bbb,\mmm)}\,\hhh,
\\{\rm or}\\
\ppp \nmid  \bbb {\rm ~~~ and~~~} (\ppp,\mmm)\mid \hhh\;.\end{array}\right.
$$
\elemm

\dem The first formula follows from the definition of
$\psi_\bbb(\ppp)$ (see Formula \eqref{eq:defpsisubd}) by considering
the three cases

$\bullet$~ $\ppp \mid \bbb$,

$\bullet$~ $\ppp \nmid \bbb$ and $\ppp\mid \mmm$, and

$\bullet$~ $\ppp\nmid \bbb$ and $\ppp\nmid \mmm$.

The second formula follows from the first one, from the definitions of
$\chi_\bbb(\ppp)$ and $\chi^*_\bbb(\ppp)$, and from the fact that
$\chi_\bbb(\ppp)\;\chi^*_\bbb(\ppp) =1$ if and only if
$\chi_\bbb(\ppp)=\chi^*_\bbb(\ppp)=1$, by considering the two cases

$\bullet$~ $\ppp \mid \bbb$ and

$\bullet$~ $\ppp \nmid \bbb$.
\cqfd

\medskip
The arithmetic function $\ccc\mapsto\mu_K(\ccc)\,\chi_\bbb(\ccc)
\;\chi^*_\bbb(\ccc) \,\Nr(\psi_\bbb(\ccc))$ being multiplicative by Lemma
\ref{lem:multiplicatives} and the complete multiplicativity of the
norm, and vanishing on the nontrivial powers of primes, the series
defining $C^*(\bbb)$ in Formula \eqref{eq:defceuler} may be written as an
Eulerian product
\begin{equation}\label{eq:prodeuler}
C^*(\bbb)=\prod_\ppp\big(1-\frac{\chi_\bbb(\ppp)
  \;\chi^*_\bbb(\ppp)\;\Nr(\psi_\bbb(\ppp))}{\Nr(\ppp)^2}\big)=
\prod_{\substack{\ppp\\\chi_\bbb(\ppp)\;\chi^*_\bbb(\ppp)=1}}
\big(1-\frac{\Nr(\psi_\bbb(\ppp))}{\Nr(\ppp)^2}\big)\;.
\end{equation}
By Equations \eqref{eq:cabktech2} and \eqref{eq:prodeuler}, and by
Lemma \ref{lem:valeurschichippsi}, we have
\[
c_{\mmm,k}=\;\frac{1}{\Nr(\mmm)}\sum_{\bbb\in \I_K^+}
\frac{\mu_K(\bbb)}{\Nr(\bbb)^2}\; \Nr((\bbb,\mmm))
\;\prod_{\ppp\,\nmid\,\bbb,\;(\ppp,\mmm)\,\mid\,\hhh}
\big(1-\frac{\Nr((\ppp, \mmm))}{\Nr(\ppp)^2}\big)
\prod_{\ppp\,\mid\, ( \bbb, \hhh)}
\big(1-\frac{1}{\Nr(\ppp)}\big)\;.
\]
Let us define $\Ga_{\mmm,\hhh}={\displaystyle
  \prod_{\substack{\ppp\\(\ppp,\mmm)\,\mid\,\hhh}}}
\big(1-\frac{\Nr((\ppp, \mmm))} {\Nr(\ppp)^2}\big)$, so that
\[
  c_{\mmm,k}=\;\frac{\Ga_{\mmm,\hhh}}{\Nr(\mmm)}\sum_{\bbb\in \I_K^+}
  \frac{\mu_K(\bbb)}{\Nr(\bbb)^2}\; \Nr((\bbb,\mmm))
  \;\prod_{\substack{\ppp\\\ppp\,\mid\,\bbb,\;(\ppp,\mmm)\,\mid\,\hhh}}
\big(1-\frac{\Nr((\ppp, \mmm))}{\Nr(\ppp)^2}\big)^{-1}
\prod_{\substack{\ppp\\\ppp\,\mid\, ( \bbb, \hhh)}}
\big(1-\frac{1}{\Nr(\ppp)}\big)\;.
\]
This equation writes $c_{\mmm,k}$ as a series
$\frac{\Ga_{\mmm,\hhh}}{\Nr(\mmm)} \sum_{\bbb\in \I_K^+}
\frac{f(\bbb)}{N(\bbb)^2}$ where $f:\I_K^+\ra\RR$ is a multiplicative
function, which vanishes on the nontrivial powers of prime ideals. By
Eulerian product, we have therefore proved Equation
\eqref{eq:simplifcmk}.

\medskip
Let us now prove that $0\leq c_{\mmm,k}\leq 1$. Note that for every
prime ideal $\ppp$, we have
\begin{equation}\label{eq:controlkappa}
1 \leq \kappa_{\mmm,\hhh}(\ppp) \leq 2 \text{ ~~and~~ }
\frac{1}{2} \leq \kappa'_\hhh (\ppp) \leq 1\;.
\end{equation}
In particular all the factors of the two products over $\ppp$ in
Equation \eqref{eq:simplifcmk} belong to $[0,1]$, hence
$0\leq c_{\mmm,k}\leq \frac{1}{\Nr(\mmm)}\leq 1$.

Let us finally prove that $c'_{\mmm}>0$. For every prime ideal $\ppp$,
let $w_\ppp= \frac{\kappa_{\mmm,\hhh}(\ppp) \; \kappa'_\hhh(\ppp)\,
  \Nr((\ppp,\mmm))}{\Nr(\ppp)^2}$. By Formula \eqref{eq:defikappa}, if
$\Nr(\ppp)=2$, we have
\[
w_\ppp=\left\{\begin{array}{l}
1/2 {\rm ~~~if~~} \ppp\mid\hhh {\rm ~~~and~~}\ppp\mid\mmm\\
1/6 {\rm ~~~if~~} \ppp\mid\hhh {\rm ~~~and~~}\ppp\nmid\mmm\\
1/2 {\rm ~~~if~~} \ppp\nmid\hhh {\rm ~~~and~~}\ppp\mid\mmm\\
1/3 {\rm ~~~if~~} \ppp\nmid\hhh {\rm ~~~and~~}\ppp\nmid\mmm
\end{array}\right.
\]
In particular $1-w_\ppp\neq 0$ if $\Nr(\ppp)=2$.
From the inequalities \eqref{eq:controlkappa} and by Equation
\eqref{eq:simplifcmk}, we have
\[
c_{\mmm,k}\geq \frac{1}{\Nr(\mmm)}
\prod_{\substack{\ppp\\(\ppp,\mmm)\,\mid\, \hhh}}
\big(1-\frac{\Nr((\ppp,\mmm))}{\Nr(\ppp)^2}\big)\;
\prod_{\ppp\,:\,\Nr(\ppp)\geq 3}
\big(1-\frac{2\,\Nr((\ppp,\mmm))}{\Nr(\ppp)^2}\big)
\prod_{\ppp\,:\,\Nr(\ppp)=2}(1-w_\ppp)\;.
\]
Since there are only finitely many primes ideals $\ppp$ dividing
$\mmm$, the term on the right hand side is bounded from below by a
positive constant $c'_\mmm=\min_{k\in\OOO_K}c_{\mmm,k}>0$. This concludes
the proof of Lemma \ref{lem:controlconstantcmmmk}. \cqfd

\medskip
Now that we understand the constant $c_{\mmm,k}\,$, we continue towards  the
proof of Theorem \ref{theo:mirsky} by giving an asymptotic formula for
the sum
\begin{equation}\label{eq:defiwtSx}
\wt S(x)=\sum_{a\in\mmm\cap C(z,\theta,x)}
\frac{\varphi_K(a)}{\Nr(a)}\;\frac{\varphi_K(a+k)}{\Nr(a+k)}\;.
\end{equation}

\blemm\label{lem:controlwtS} Uniformly in $\mmm\in \I_K^+$,
$k\in\OOO_K$, $z\in \CC^\times$, $\theta\in \;]0,2\pi]$ and $x\geq 1$,
we have
\begin{equation}\label{eq:asymwtS}
  \wt S(x)=\frac{\theta\;c_{\mmm,k}}{\sqrt{|D_K}}\;x^2+
  \bigO(x)\;.
\end{equation}
\elemm

\dem For all nonzero elements $a$ and $b$ in the factorial ring
$\OOO_K$, we denote by $(a,b)$ any fixed choice of gcd of $a$ and $b$,
and by $[a,b]$ any fixed choice of lcm of $a$ and $b$.

By Equation \eqref{eq:MobiusinvEuler}, for every $a\in\OOO_K-\{0\}$,
we have
\[
\frac{\varphi_K(a)}{\Nr (a)}=\frac{1}{|\OOO_K^\times|}
\sum_{b\in\OOO_K-\{0\}\;:\;b\,\mid\,a}\;\frac{\mu_K(b)}{\Nr (b)}\;.
\]
Let $x\geq 1$. Applying twice this equality, since $\Nr(b)\leq \Nr(a)$
when $b\mid a$, we have by Fubini's theorem
\begin{align}
  \wt S(x)&=\frac{1}{|\OOO_K^\times|^2}\sum_{a\in\mmm\cap C(z,\theta,x)}
  \;\;\sum_{b\in\OOO_K-\{0\}\,:\;b\,\mid\,a}\;\frac{\mu_K(b)}{\Nr (b)}
  \sum_{c\in\OOO_K-\{0\}\,:\;c\,\mid\,a+k}\;\frac{\mu_K(c)}{\Nr (c)}
  \nonumber\\&=\frac{1}{|\OOO_K^\times|^2}
  \sum_{b\in\OOO_K-\{0\}\,:\;|b|\leq x}\;\frac{\mu_K(b)}{\Nr (b)}
  \sum_{c\in\OOO_K-\{0\}}\;\frac{\mu_K(c)}{\Nr (c)}\sum_{\substack{a
      \in\mmm\cap C(z,\theta,x)\\b\,\mid\,a,\;c\,\mid\,a+k}} 1\;.
  \label{eq:firstreducwtS}
\end{align}
Let $b,c\in\OOO_K-\{0\}$. The system of three congruences
${\displaystyle \left\{\begin{array}{l} a\equiv0\!\!\mod
  m\\ a\equiv0\!\!\mod b\\a\equiv -k\!\!\mod c\end{array}\right.}$ has
a solution $a\in\OOO_K-\{0\}$ such that $|a|\leq x$ if and only if
there exists an element $n\in\OOO_K-\{0\}$ such that $a=bn$, $|n|\leq
\frac{x}{|b|}$ and
\begin{equation}\label{eq:systdeucongru}
\left\{\begin{array}{l} bn\equiv 0\!\!\!\mod m\\ bn\equiv -k\!\!\!\mod
c\;.\end{array} \right.
\end{equation}
When $(b,c)\nmid k$, no solution exists.

Assume that $(b,c)\mid k$. Since $\frac{b}{(b,c)}$ is invertible
modulo $\frac{c}{(b,c)}$, we denote by $\overline{\frac{b}{(b,c)}}$ a
multiplicative inverse of $\frac{b}{(b,c)}$ modulo $\frac{c}{(b,c)}$.
Then the system of congruences \eqref{eq:systdeucongru} is equivalent to
\begin{equation}\label{eq:bicong}
  \left\{\begin{array}{l} \frac{b}{(b,m)}n\equiv 0\!\!\!
  \mod \frac{m}{(b,m)}\\ \frac{b}{(b,c)}n\equiv -\frac{k}{(b,c)}\!\!\!\mod
  \frac{c}{(b,c)}\end{array} \right. \Leftrightarrow\;
  \left\{\begin{array}{l} n\equiv 0\!\!\!\mod \frac{m}{(b,m)}\\
  n\equiv -\frac{k}{(b,c)}\,\overline{\frac{b}{(b,c)}}\!\!\!\mod
\frac{c}{(b,c)}\;.\end{array} \right.
\end{equation}

Recall that a system of two congruences $\left\{\begin{array}{l}
n\equiv \alpha_0\!\!\mod \alpha\\ n\equiv \beta_0\!\!\mod
\beta \end{array} \right.$ with unknown $n\in\OOO_K$, where $\alpha,
\beta,\alpha_0,\beta_0\in \OOO_K$ and $\alpha,\beta\neq 0$, has a
solution if and only if $\alpha_0-\beta_0\equiv 0\!\!  \mod(\alpha,
\beta)$.  Furthermore, if this congruence condition is satisfied, that
is, if there exists $n_0,m_0\in \OOO_K$ such that $\alpha_0-\beta_0=
\beta m_0-\alpha n_0$, then $n$ is a solution if and only if
\[
n-\alpha_0-\alpha n_0\in\alpha\OOO_K\cap\beta\OOO_K
=[\alpha,\beta]\OOO_K\;.
\]
This is equivalent to asking $n$ to belong to the translate
$\Lambda_{\alpha,\beta,\alpha_0,\beta_0}= \alpha_0+ \alpha n_0 +
\vec\Lambda_{\alpha,\beta}$ of the $\ZZ$-lattice
$\vec\Lambda_{\alpha,\beta}=[\alpha,\beta]\OOO_K$.

Applying this with $\alpha=\frac{m}{(b,m)}$, $\beta=\frac{c}{(b,c)}$,
$\alpha_0=0$ and $\beta_0= -\frac{k}{(b,c)}\, \overline{\frac{b}
  {(b,c)}}$, since the elements $\frac{b}{(b,c)}$ and
$\frac{b}{(b,m)}$ are both coprime with $\big( \frac{m} {(b,m)},
\frac{c}{(b,c)}\big)$, the system \eqref{eq:bicong} has a solution if
and only if the following divisibility condition holds
\begin{align*}
&\Big(\frac{m}{(b,m)},\frac{c}{(b,c)}\Big)\; \mid
\frac{k}{(b,c)}\;\overline{\frac{b}{(b,c)}}\Leftrightarrow
\Big(\frac{m}{(b,m)},\frac{c}{(b,c)}\Big)\; \mid
\frac{k}{(b,c)}\\ \Leftrightarrow\;\;&
\Big(\frac{m}{(b,m)},\frac{c}{(b,c)}\Big)\; \mid
\frac{k}{(b,c)}\;\frac{b}{(b,m)}\Leftrightarrow
\big(m(b,c),c(b,m)\big)\, \mid
k\,b\;.
\end{align*}
Thus Equation \eqref{eq:firstreducwtS} becomes, using Equation
\eqref{eq:dilatCzthetaR},
\[%\begin{align}
  \wt S(x)=\frac{1}{|\OOO_K^\times|^2}
  \sum_{\substack{b,c\,\in\OOO_K-\{0\}\,:\;|b|\leq x\\
      (b,c)\,\mid\, k,\;(m(b,c),c(b,m))\, \mid \,k\,b}}\;
  \frac{\mu_K(b)\,\mu_K(c)}{\Nr (b)\,\Nr (c)}
\sum_{n\in\Lambda_{\alpha,\beta,\alpha_0,\beta_0}\cap C(b^{-1}z,\,\theta,\,x/|b|)} 1\;.
\]%  \label{eq:secondreducwtS}\end{align}
Let $b,c$ be as in the index of the first sum above. Using again the
standard Gauss counting argument, using Formula
\eqref{eq:calccocoldiamideal} for the second equality and the equation
$\Nr([\alpha,\beta])=\frac{\Nr(\alpha)\Nr(\beta)}
      {\Nr((\alpha,\beta))}$ for the last equality, we have, uniformly
      in $b,c,m\in\OOO_K-\{0\}$, $k\in\OOO_K$, $z\in\CC^\times$,
      $\theta\in\;]0,2\pi]$ and $y\geq 1$,
\begin{align*}
  \card&(\Lambda_{\alpha,\beta,\alpha_0,\beta_0}\cap
  C(b^{-1}z,\theta,y))=
  \frac{\theta}{2\,\covol_{\vec\Lambda_{\alpha,\beta}}}\;y^2
  +\bigO\Big(\;\frac{\diam_{\vec\Lambda_{\alpha,\beta}}}
  {\covol_{\vec\Lambda_{\alpha,\beta}}}\;y\Big)\\  =\;\;&
  \frac{\theta}{\sqrt{|D_K|}\Nr([\frac{m}{(b,m)},
      \frac{c}{(b,c)}])}\;y^2  +\bigO\Big(\;\frac{1}
  {\sqrt{\Nr([\frac{m}{(b,m)},\frac{c}{(b,c)}])}}\;y\Big)\\=\;\;&
  \frac{\theta\;\Nr\!\big((m(b,c),c(b,m))\big)}
       {\sqrt{|D_K|}\,\Nr(m)\,\Nr(c)}\;y^2  +
\bigO\Big(\;\frac{\Nr\!\big((m(b,c),c(b,m))\big)^{1/2}}
  {\Nr(m)^{1/2}\,\Nr(c)^{1/2}}\;y\Big) \;.
\end{align*}
Using this with $y=\frac{x}{|b|}$, which is at least $1$
since $|b|\leq x$, we have
\begin{align}
\wt S(x)&=\frac{\theta\,x^2}{\sqrt{|D_K|}}
\sum_{\substack{b,c\,\in\OOO_K-\{0\}\,:\;|b|\leq x\\
(b,c)\,\mid\, k,\;(m(b,c),c(b,m))\, \mid \,k\,b}}\;
\frac{\mu_K(b)\,\mu_K(c)\,\Nr\!\big((m(b,c),c(b,m))\big)}
{|\OOO_K^\times|^2\,\Nr (b)^2\,\Nr (c)^2\,\Nr(m)}\nonumber
\\ &\;\;\;+\bigO\Big(x\sum_{b,c\,\in\OOO_K-\{0\}}\;
\frac{\Nr\!\big((m(b,c),c(b,m))\big)^{1/2}}
{|\OOO_K^\times|^2\,\Nr (b)^{3/2}\,\Nr (c)^{3/2}\,\Nr(m)^{1/2}}
\Big)\;.\label{eq:secondreducwtS}
\end{align}

By Equation \eqref{eq:Zmmmk} (replacing therein $x$ by $x^2$),
completing the first sum of the above equation with the indices
$b\in\OOO_K-\{0\}$ such that $|b|> x$ introduces an error of the form
$\bigO(\frac{1}{x})$ (uniformly in $m\in\OOO_K-\{0\}$, $k\in\OOO_K$
and $x\geq 1$). A computation similar to the one done for Equation
\eqref{eq:Zmmmk} gives that the second sum in Equation
\eqref{eq:secondreducwtS} is actually bounded by
$\frac{1}{|\OOO_K^\times|^2} \zeta_K(\frac{3}{2})^2\,\zeta_K(2)$, which
is uniform since there are only finitely many such fields $K$.

By the definition of the constant $c_{\mmm,k}$ in Equation
\eqref{eq:deficmk}, this proves Equation \eqref{eq:asymwtS}, hence
concludes the proof of Lemma \ref{lem:controlwtS}.
\cqfd

\medskip
\noindent{\bf Proof of Theorem
  \ref{theo:mirsky}.} For all $a,k\in\OOO_K$ with $a\neq 0$, we have
\[
\Nr(a+k)=\Nr(a)\Big|1+\frac{k}{a}\Big|^2\leq \Nr(a)
\Big(1+2\sqrt{\frac{\Nr(k)}{\Nr(a)}}+\frac{\Nr(k)}{\Nr(a)}\Big)\;,
\]
and similarly $\Nr(a+k)\geq \Nr(a) \big(1-2\sqrt{\frac{\Nr(k)}
  {\Nr(a)}}+\frac{\Nr(k)}{\Nr(a)}\big)$. Let us define the maps
$f_\pm:[1,+\infty[\;\ra\RR$ by $t\mapsto t^2\pm 2\sqrt{\Nr(k)}
    \,t^{3/2}+\Nr(k)\,t$, so that their derivatives are $f'_\pm(t)=2t\pm
    3\sqrt{\Nr(k)} \,t^{1/2}+\Nr(k)$ and
\begin{equation}\label{eq:proprifpm}
\frac{f_-(\Nr(a))}{\Nr(a)\,\Nr(a+k)}\leq 1\leq
\frac{f_+(\Nr(a))}{\Nr(a)\,\Nr(a+k)}
\end{equation}
    
For all $z\in\CC^\times$, $\theta\in\;]0,2\pi]$, $x\geq 1$ and
$n\in\NN-\{0\}$, let
\[
a_n=\sum_{a\in\mmm\cap C(z,\theta,x\,)\,:\,\Nr(a)=n}
\frac{\varphi_K(a)}{\Nr(a)}\;\frac{\varphi_K(a+k)}{\Nr(a+k)}\;,
\]
so that by Equation \eqref{eq:defiwtSx}, we have ${\displaystyle \wt
  S(x)=\sum_{1\leq n\leq x^2} a_n}$.

By the definition \eqref{eq:defiSzthetakmmmx} of the sum
$S_{z,\theta,k,\mmm}(x)$ and the inequalities \eqref{eq:proprifpm}, by
Abel's summation formula, by applying twice Lemma
\ref{lem:controlwtS}, and since $c_{\mmm,k}\leq 1$ by Lemma
\ref{lem:controlconstantcmmmk}, we have
\begin{align*}
&S_{z,\theta,k,\mmm}(x)\leq \sum_{1\leq n\leq x^2}
a_n\,f_+(n)=\Big(\sum_{1\leq n\leq x^2} a_n\Big)\;f_+(x^2) -
\int_{1}^{x^2}\Big(\sum_{1\leq n\leq t} a_n\Big)\;f'_+(t)\;dt\\=\; &
\Big(\frac{\theta\;c_{\mmm,k}}{\sqrt{|D_K|}}\;x^2+
\bigO(x)\Big)\Big(x^4+ 2\sqrt{\Nr(k)}\,x^{3}+\Nr(k)\,x^2\Big)
\\ &- \int_{1}^{x^2}\Big(\frac{\theta\;c_{\mmm,k}}{\sqrt{|D_K|}}\;t+
  \bigO(t^{1/2})\Big)\Big(2t+ 3\sqrt{\Nr(k)}
  \,t^{1/2}+\Nr(k)\Big)\;dt\\=\; &
  \frac{\theta\;c_{\mmm,k}}{3\,\sqrt{|D_K|}}\;x^6+
  \bigO\big((1+\sqrt{\Nr(k)}\,)\,x^5+\Nr(k)\,x^4\big)\;.
\end{align*}
Replacing $f_+$ by $f_-$ gives the same minoration to
$S_{z,\theta,k,\mmm}(x)$, hence Theorem \ref{theo:mirsky} follows.
\cqfd

{\small \bibliography{../biblio} }
%{\small \bibliography{viitteet} }

\bigskip
{\small
\noindent \begin{tabular}{l} 
Department of Mathematics and Statistics, P.O. Box 35\\ 
40014 University of Jyv\"askyl\"a, FINLAND.\\
{\it e-mail: jouni.t.parkkonen@jyu.fi}
\end{tabular}
\medskip

\noindent \begin{tabular}{l}
Laboratoire de mathématique d'Orsay, UMR 8628 CNRS,\\
Universit\'e Paris-Saclay,\\
91405 ORSAY Cedex, FRANCE\\
{\it e-mail: frederic.paulin@universite-paris-saclay.fr}
\end{tabular}}

\end{document}